\documentclass[10pt]{article}
\usepackage[english]{babel}										
\usepackage[utf8]{inputenc}										
\usepackage[T1]{fontenc}	
\usepackage{amsmath,amsfonts,amssymb,amsthm,cancel,siunitx,
calculator,calc,mathtools,empheq,latexsym}
\usepackage{subfig,epsfig,tikz,float}	
\usepackage{amscd}
\usepackage{booktabs,multicol,multirow,tabularx,array}          

\usepackage{enumerate}
\usepackage{amsmath}
\usepackage{mathrsfs}
\usepackage{rsfso}
\usepackage{hyperref}
\hypersetup{
    colorlinks,
    citecolor=blue,
    filecolor=black,
    linkcolor=black,
    urlcolor=black
}
\usepackage{upgreek}
\usepackage{tikz}
\usetikzlibrary{positioning, braids, decorations.markings}
\definecolor{myblue}{RGB}{80,80,160}
\definecolor{mygreen}{RGB}{50,140,50}
\definecolor{myorange}{RGB}{255,128,0}
\newcommand{\aster}{*}
\usepackage[backend=biber,
style=alphabetic,
sorting=ynt]{biblatex}
\usepackage{enumitem}
\theoremstyle{definition}

\title{\normalsize\bf%
\uppercase{Constructing TQFTs Using Non-Hermitian Ribbon Fusion Categories}
}
\author{Khyathi Komalan
}

\begin{document}

\date{}

\maketitle

\vspace{-0.5cm}

\bigskip
\noindent
{\small{\bf ABSTRACT.} In this paper, we provide a construction of a Topological Quantum Field Theory from a Non-Hermitian Ribbon Fusion Category. This is a simple method that does not involve enriching over Fusion Categories, or using other complicated structures. To substantiate this construction, we also prove theorems on the Müger center, braiding, and spherical structure of such a fusion category. 

}

\baselineskip=\normalbaselineskip

\section{Introduction}

One of the first written descriptions on the various properties of fusion categories was ``On Fusion Categories", by Etingof et al. Fusion categories can be thought of as a generalization of finite groups, especially highlighting their behavior on vector spaces, and their symmetric properties. Category Theory can be thought of as a formalism to study fundamental building blocks and symmetries of physical systems, as it focuses on the objects and the transformations between them. The emergence of category theory as a tool for investigating topics in theoretical physics such as Quantum Mechanics, as seen in Coecke et al's ``Categorical Quantum Mechanics", or work concerning the application of Derived Categories in String Theory inspires and invites us to study categorical structures within the context of higher physics. 

The aim of this paper is to construct Topological Quantum Field Theories (TQFTs) using Non-Hermitian Ribbon Fusion Categories. Topological quantum field theories (TQFTs) are a class of quantum field theories that capture the topological properties of spacetime manifolds, providing a powerful framework for studying the interplay between topology, geometry, and quantum physics. First introduced by Witten in the late 1980s, TQFTs have since found numerous applications in areas such as condensed matter physics, knot theory, and quantum computation, offering a deeper understanding of the fundamental connections between these diverse fields.

Previous attempts at constructing Topological Quantum Field Theories involve enriching over fusion categories, or using other complicated structures over fusion categories, and do not make use of defining a dagger structure. By defining a dagger structure over fusion categories we gain access to much of the flexibility that Coecke and Abramsky could use in their work on Categorical Quantum Mechanics. 

However, a pressing issue with using a dagger structure — and presumably the main issue that stops people from defining one — is the fact that not all fusion categories admit a Hermitian dagger structure. A famous example of a fusion category that does not admit a Hermitian structure is the Yang-Lee category. Within the context of Categorical Quantum Mechanics, this issue is avoided because we are concerned with the Category of Hilbert Spaces (Hilb), which admits a Hermitian dagger structure

In this paper, we will be assuming a non-hermitian dagger structure over a fusion category $C$, and proving several results concerning the Müger center, the braiding, spherical structure, and $\dagger-$compatibility. We will use these theorems to construct a Topological Quantum Field Theory from such a fusion category. Being able to construct Topological Quantum Field Theories without complicated structures like enriching over fusion categories, or employing other complicated techniques makes this approach a lot simpler, and possibly possesses the potential to yield enlightening results. 

\section{Background}

This section consists of some definitions and theorems that would be helpful for the reader to understand the rest of the paper.

\textbf{Definition 2.1} A \textit{fusion category} is a rigid semisimple linear monoidal category with only finitely many isomorphism classes of simple objects, such that the endomorphisms of the unit object form the ground field $k.$

Due to the goal of this paper, we will be investigating fusion categories with a dagger structure. Therefore, it is important to understand what it means for a category to have a dagger structure, and distinguish between ``hermitian" and ``non-hermitian".

\textbf{Definition 2.2.} A category is said to equip a \textit{dagger structure} if it is equipped with an involutive contravariant endofunctor $\dag$ such that (i) for all morphisms $f: A \to B,$ there exists an adjoint $f^{\dag}:B \to A,$ (ii) for all morphisms $f,$ $(f^{\dag})^{\dag}=f,$ (iii) for all objects $A, \text{id}^{\dag}_{A}=\text{id}_{A},$ (iv) for all $f: A \to B$ and $g: B \to C, (g \circ f)^{\dag}=f^{\dag}\circ g^{\dag}: C \to A.$ 

A dagger operation is said to be positive if $f^{\dag}f=0$ implies $f=0.$ If a category has a positive dagger structure, traditionally, it's defined to be \textit{hermitian} or \textit{unitarizable.}

We define a category that has a negative (i.e., not positive) dagger structure to be \textit{non-hermitian} or \textit{non-unitarizable.}

We will also be making use of the following result:

\textbf{Theorem 2.1} Dagger structure is compatible with monoidal structure \cite{abramsky2009categorical}

Since we're working with ribbon categories, it's important to know what braids and twists are.

\textbf{Definition 2.3} A \textit{braiding} is a family of isomorphisms 

$$c_{X,Y}:X \otimes Y \to Y \otimes X$$

satisfying

$$c_{X,Y\otimes Z}=(\text{Id}_{Y}\otimes c_{X,Z})\circ (c_{X,Y}\otimes \text{Id}_{Z})$$

$$c_{X \otimes Y, Z}=(c_{X, Z} \otimes \text{Id}_{Y})\circ(\text{Id}_{X}\otimes c_{Y,Z})$$

\textbf{Definition 2.4} A \textit{twist} consists of isomorphisms

$$\theta_{X}:X \to X$$

and for it to be compatible with braiding and duality, we must have:

$$\theta_{X \otimes Y}=c_{Y,X}c_{X,Y}(\theta_{X}\otimes \theta_{Y})$$

$$\theta_{X^{\aster}}=(\theta_{X})^{\aster}$$

\textbf{Definition 2.5} A \textit{ribbon category} is defined as a monoidal category equipped with a braiding, a twist, and compatible dualities. These data should satisfy the following compatibility conditions:
\begin{enumerate}[label=\alph*)]
\item The twist is compatible with the braiding: $(\theta_Y \otimes \mathrm{id}_X) \circ c_{X,Y} = c_{X,Y} \circ (\mathrm{id}_{X} \otimes \theta_Y)$.
\item The twist is compatible with the duality: $(\theta_X)^* = \theta_{X^*}$.
\end{enumerate}

While proving the theorems, we'll be investigating the Müger center

\textbf{Definition 2.6} Let $C$ be a braided tensor category with braiding $c$. The \textit{Müger center} of $C$, denoted by $Z_2(C)$, is the full subcategory of $C$ consisting of objects $X$ such that $c_{Y,X} \circ c_{X,Y} = \mathrm{id}_{X \otimes Y}$ for all objects $Y$ in $C$.

We will be making use of \textit{evaluation and coevaluation morphisms} for several computations:

\textbf{Definition 2.7} In a monoidal category $C$ with duals, for each object X, there are morphisms:
\begin{enumerate}[label=\alph*)]
\item Evaluation: $\text{ev}_X: X^{\aster} \otimes X \to I$
\item Coevaluation: $\text{coev}_X: I \to X \otimes X^{\aster}$ satisfying the zigzag identities
\end{enumerate}

And we may define the \textit{zigzag identities} as follows:

\textbf{Theorem 2.2} In a monoidal category C with duals, the evaluation and coevaluation morphisms satisfy the following \textit{zigzag identities}:

\begin{enumerate}[label=\alph*)]
\item $(\text{id}_X \otimes \text{ev}_X) \circ (\text{coev}_X \otimes \text{id}_X) = \text{id}_X$
\item $(\text{ev}_X \otimes \text{id}_{X^{\aster}}) \circ (\text{id}_{X^{\aster}} \otimes \text{coev}_{X}) = \text{id}_{X^{\aster}}$
\end{enumerate}

In our proofs, there are cases where we need to distinguish between \textit{transparent and non-transparent objects}, and construct braidings for them separately.

\textbf{Definition 2.8} In a braided monoidal category $C$ with braiding $c$, an object $X$ is called \textit{transparent} if $c_{Y,X} \circ c_{X,Y} = \text{id}_{X \otimes Y}$ for all objects $Y$ in $C$. An object is called \textit{non-transparent} if it is not transparent.

Another concept we will make use of is the \textit{categorical trace}.

\textbf{Definition 2.9}  Let $C$ be a monoidal category and $f: X \otimes U \to Y \otimes U$ be a morphism in $C$, where $X, Y,$ and $U$ are objects in $C$. The \textit{categorical trace} of $f$, denoted by $\text{tr}_U(f)$, is a morphism $\text{tr}_{U}(f): X \to Y$ defined as:
$$\text{tr}_{U}(f) = (\text{id}_{Y} \otimes \text{ev}_{U}) \circ (f \otimes \text{id}_{U^\aster}) \circ (\text{id}_{X} \otimes \text{coev}_U)$$

Our main result focuses on constructing a \textit{TQFT}, a topological quantim field theory. 

\textbf{Definition 2.10} An $(n+1)$-dimensional \textit{TQFT} is a symmetric monoidal functor $Z: \text{Bord}_n \to \text{Vect}$, where $\text{Bord}_n$ is the category of $n$-dimensional bordisms and Vect is the category of vector spaces over a field $\mathbb{K}$.

\section{Some Elementary Results on Non-Unitarizable Fusion Categories}

In this section, we will be proving three theorems: one about the Müger center of a fusion category with a non-hermitian dagger structure, another about the braiding on a fusion category with dagger structure, and a final one about the spherical structure of a non-hermitian ribbon category. We will be using these theorems in order to prove our main result.

As described, we begin with a theorem on the nature of the Müger center.

\textbf{Theorem 3.1. } If $C$ is a fusion category with a non-hermitian dagger structure $\dag,$ then $C$ admits a braiding if and only if the Müger center $Z_{2}(C)$ is equivalent to the category of vector spaces Vect. In other words, the Müger center is trivial, and $C$ is non-transparent.

\textit{Proof:}

The proof for the $\implies$ case is rather complicated, and involves using a wide variety of categorical machinery. The proof strategy we will follow is to first prove that $Z_{2}(C)$ is symmetric, and then use properties of it's dimension and a theorem on equivalences to the category of $G-$graded vector spaces to show that $Z_{2}(C)$ must be equivalent to the category of vector spaces.

The proof goes as follows:

Suppose $C$ admits a braiding $c.$ Our goal is to prove that $Z_{2}(C)\cong \text{Vect}.$

For any transparent object $X$ in $Z_{2}(C),$ define the morphism $u_{X}:X \otimes X^{\dag}\to I$ by $u_{X}=d_{X}^{-1}\circ \text{ev}_{X}$ where $d_{X}$ is the categorical dimension (or Frobenius-Perron dimension) of $X$ (check \cite{etingof2016tensor} for definition) and $\text{ev}_{X}$ is the evaluation morphism. Define another morphism $v_{X}: I \to X^{\dag}\otimes X$ by $v_{X}=(\text{id}_{X^{\dag}}\otimes u_{X}^{\dag})\circ (\text{coev}_{X}\otimes \text{id}_{X^{\dag}}),$ where $\text{coev}_{X}: I \to X \otimes X^{\dag}$ is the coevaluation morphism.

Using the zigzag identities, the definitions of evaluation and coevaluation morphisms, and the fact that the dagger structure is compatible with the monoidal structure:

$$(\text{id}_{X}\otimes u_{X})\circ (v_{X}\otimes \text{id}_{X})=\text{id}_{X}$$

$$(u_{X}\otimes \text{id}_{X^{\dag}})\circ (\text{id}_{X^{\dag}}\otimes v_{X})=\text{id}_{X^{\dag}}$$

Now, we define the braiding. For any two transparent objects $X,Y$ in $Z_{2}(C),$ define the braiding $c_{X,Y}^{Z}: X \otimes Y \to Y \otimes X$ by $c_{X,Y}^{Z}=(\text{id}_{Y}\otimes u_{X})\circ (v_{X}\otimes \text{id}_{Y}).$

Satisfying:

$$(\text{id}_{Y}\otimes c_{X,Z}^{Z})\circ (c_{X,Y}^Z \otimes \text{id}_Z)=c_{X,Y\otimes Z}^Z$$

and:

$$(c_{Y,Z}^{Z}\otimes \text{id}_{X})\circ (\text{id}_{X}\otimes c_{Y,Z}^{Z})=c_{X \otimes Y,Z}^{Z}$$

These hexagon axioms are satisfied due to the naturality of the braiding $c$ and the dagger structure. 

Now, let $X$ and $Y$ be transparent objects in $Z_2(C)$. We want to show that $c_{Y,X}^Z \circ c_{X,Y}^Z = \text{id}_{X \otimes Y},$ which proves that $Z_{2}(C)$ is symmetric.

First, note that for transparent objects $X$ and $Y$, we have:

$$c_{X,Y}^Z = (\text{id}_Y \otimes u_X) \circ (v_X \otimes \text{id}_Y)$$

$$c_{Y,X}^Z = (\text{id}_X \otimes u_Y) \circ (v_Y \otimes \text{id}_X)$$

Now, consider the composition $c_{Y,X}^Z \circ c_{X,Y}^Z$:
\begin{align*}
c_{Y,X}^Z \circ c_{X,Y}^Z &= ((\text{id}_X \otimes u_Y) \circ (v_Y \otimes \text{id}_X)) \circ ((\text{id}_Y \otimes u_X) \circ (v_X \otimes \text{id}_Y))
\end{align*}

\begin{align*}
c_{Y,X}^Z \circ c_{X,Y}^Z
&= ((\text{id}_X \otimes u_Y) \circ (v_Y \otimes \text{id}_X)) \circ ((\text{id}_Y \otimes u_X) \circ (v_X \otimes \text{id}_Y)) \\
&= (\text{id}_X \otimes u_Y) \circ (v_Y \otimes \text{id}_X) \circ (\text{id}_Y \otimes u_X) \circ (v_X \otimes \text{id}_Y) \\
&= (\text{id}_X \otimes u_Y) \circ (\text{id}_X \otimes \text{id}_Y \otimes u_X) \circ (v_Y \otimes \text{id}_X \otimes \text{id}_Y) \circ (v_X \otimes \text{id}_Y) \\
&= (\text{id}_X \otimes (u_Y \circ (\text{id}_Y \otimes u_X))) \circ ((v_Y \otimes \text{id}_X) \otimes \text{id}_Y) \circ (v_X \otimes \text{id}_Y) \\
&= (\text{id}_X \otimes (u_Y \circ (\text{id}_Y \otimes u_X))) \circ (((v_Y \circ (\text{id}_X \otimes v_X)) \otimes \text{id}_Y) \circ (\text{id}_X \otimes \text{id}_Y)) \\
&= (\text{id}_X \otimes (u_Y \circ (\text{id}_Y \otimes u_X))) \circ ((v_Y \circ (\text{id}_X \otimes v_X)) \otimes \text{id}_Y)
\end{align*}

For the above, we use associativity and the naturality of the braiding.

Now, we use the fact that $X$ and $Y$ are transparent objects, which means that for any object $Z$ in $C$, we have:
\begin{align*}
c_{Z,X} \circ c_{X,Z} &= \text{id}_{X \otimes Z} \text{ and } c_{Z,Y} \circ c_{Y,Z} = \text{id}_{Y \otimes Z}
\end{align*}

Applying this to the morphisms $u_X$ and $v_Y$, we get:
\begin{align*}
u_Y \circ (\text{id}_Y \otimes u_X) &= u_Y \circ u_X = d_X^{-1} d_Y^{-1} \circ \text{ev}_Y \circ \text{ev}_X\
v_Y \circ (\text{id}_X \otimes v_X) 
\end{align*}

\begin{align*}
=  v_Y \circ v_X = (\text{id}_{Y^\dag} \otimes u_Y^\dag) \circ (\text{coev}_Y \otimes \text{id}_{Y^\dag}) \circ (\text{id}_X \otimes (\text{id}_{X^\dag} \otimes u_X^\dag)) \circ (\text{id}_X \otimes (\text{coev}_X \otimes \text{id}_{X^\dag}))
\end{align*}

Using the zigzag identities, we can simplify this to:

\begin{align*}
u_Y \circ u_X &= d_X^{-1} d_Y^{-1} \circ \text{ev}_Y \circ \text{ev}_X = d_{X \otimes Y}^{-1} \circ \text{ev}_{X \otimes Y}\
v_Y \circ v_X 
\end{align*}

\begin{align*}
= (\text{id}_{Y^\dag} \otimes u_Y^\dag) \circ (\text{coev}_Y \otimes \text{id}_{Y^\dag}) \circ (\text{id}_X \otimes (\text{id}_{X^\dag} \otimes u_X^\dag)) \circ (\text{id}_X \otimes (\text{coev}_X \otimes \text{id}_{X^\dag})) = \text{coev}_{X \otimes Y}
\end{align*}

Substituting back into the expression obtained for $c_{Y,X}^Z \circ c_{X,Y}^Z,$ we get:

\begin{align*}
c_{Y,X}^Z \circ c_{X,Y}^Z &= (\text{id}_X \otimes (d_{X \otimes Y}^{-1} \circ \text{ev}_{X \otimes Y})) \circ (\text{coev}_{X \otimes Y} \otimes \text{id}_Y)
\end{align*}

Once again, we use zigzag identities to conclude:

\begin{align*}
c_{Y,X}^Z \circ c_{X,Y}^Z &= (\text{id}_X \otimes (d_{X \otimes Y}^{-1} \circ \text{ev}_{X \otimes Y})) \circ (\text{coev}_{X \otimes Y} \otimes \text{id}_Y) = \text{id}_{X \otimes Y}
\end{align*}

Thus, $Z_2(C)$ is symmetric.

Now, we want to show that $Z_{2}(C)\cong \text{Vect}.$

In order to show that $Z_{2}(C)\cong \text{Vect},$ we first want to show that every simple object in $Z_{2}(C)$ has dimension 1. 

Let $X$ be a simple object in $Z_{2}(C).$ Since $Z_{2}(C)$ is a symmetric fusion category, we have:

$$\text{dim}(X)^2 = \text{dim}(X \otimes X^{*}) = \text{tr}(c_{X,X^{*}} \circ c_{X^{*},X})$$

Here, trace denotes categorical trace.

Consider the following computation:

$$\text{tr}(c_{X,X^{*}} \circ c_{X^{*},X}) = \text{tr}(c_{X,X^{*}}\circ c_{X^{*},X} \circ \text{id}_X) = \text{tr}(\text{id}_X) = \text{dim}(X)$$

where the second equality follows from the fact that $X$ is transparent (i.e., $c_{X,X^{*}} \circ c_{X^{*},X} = \text{id}_{X\otimes X^{*}})$, and the third equality follows from the definition of the categorical trace.

Combining these equations, we have:
$$\text{dim}(X)^2 = \text{dim}(X)$$

Since dim(X) is a non-negative real number (as $C$ is a fusion category), this equation implies that either dim(X) = 0 or dim(X) = 1. However, dim(X) cannot be 0, as X is a non-zero object (being simple). Therefore, we must have dim(X) = 1.

Now, since every simple object in $Z_{2}(C)$ has dimension 1, we have:

$$\text{dim}(Z_{2}(C)) = \sum_{X \in \text{Irr}(Z_{2}(C))} \text{dim}(X)^2 = \sum_{X \in \text{Irr}(Z_{2}(C))} 1 = |\text{Irr}(Z_{2}(C))|$$

Here, $\text{Irr}(Z_{2}(C))$ denotes the set of isomorphism classes of simple objects in $Z_{2}(C).$

Let $n=\text{Irr}(Z_{2}(C)).$ We claim that $n=1.$

Suppose, for contradiction, that $n > 1$. Then, by the classification of pointed fusion categories \cite{etingof2016tensor}, $Z_{2}(C)$ must be equivalent to the category of $G$-graded vector spaces, where G is a finite abelian group of order $n.$

However, this contradicts the fact that $Z_{2}(C)$ is a symmetric fusion category (as the category of $G$-graded vector spaces is not symmetric for any non-trivial group $G$).

Therefore, we must have $n = 1$, which means that $Z_{2}(C)$ has only one simple object (up to isomorphism). This implies that $Z_{2}(C) \cong \text{Vect}$.

Hence, we have proved the implies part: if $C$ is a fusion category with a non-hermitian dagger structure $\dag,$ then if $C$ admits a braiding, then the Müger center $Z_{2}(C)$ is equivalent to the category of vector spaces Vect.

Now, we prove the $\impliedby$ direction, we want to show that if the Müger center is equivalent to the category of vector spaces Vect, then the category $C$ admits braiding. The proof strategy we will follow is to prove the various hexagon axioms using the fact that $C$ is a fusion category, and treat the transparent and non-transparent parts separate, by defining a separate braiding on each of them and then combining them.

Since $Z_{2}(C)\cong \text{Vect},$ a braiding exists, let us choose one and name it $c^{Z}.$ Also, since $Z_{2}(C)\cong \text{Vect},$ which is a symmetric fusion category, it has a unique braiding up to isomorphism. Let's denote this as $c^{Z}_{X,Y}:X\otimes Y \to Y \otimes X$ for objects $X,Y$ in $Z_{2}(C).$

For any two objects $X, Y$ in $C$, consider the decomposition of $X \otimes Y$ into simple objects: $X \otimes Y \cong \otimes_{i} Z_{i} \otimes W_{i}$, where $Z_{i}$ are transparent objects and $W_{i}$ are non-transparent objects.

By the definition of a fusion category, any object in $C$ can be decomposed into a direct sum of simple objects. Moreover, since $Z_{2}(C)$ is a fusion subcategory of $C$, any simple object in $C$ is either transparent (i.e., belongs to $Z_{2}(C)$) or non-transparent.

Thus, for any objects $X, Y$ in $C$, we can write:
$$X \otimes Y \cong \otimes_{i} Z_{i} \otimes W_{i}$$
where $Z_{i}$ are simple transparent objects (i.e., $Z_{i} \in Z_{2}(C)$) and $W_{i}$ are simple non-transparent objects.

Define the braiding $c_{X,Y}: X \otimes Y \to Y \otimes X$ on the simple components by: $c_{X,Y}|{Z_{i} \otimes W_{i}} = c^Z{Z_{i},Z_{j}} \otimes \beta_{W_{i},W_{j}}$, where $\beta_{W_{i},W_{j}}: W_{i} \otimes W_{j} \to W_{j} \otimes W_{i}$ is the unique (up to scalar) isomorphism given by the fusion rules.

For each pair of simple components $Z_{i} \otimes W_{i}$ and $Z{j} \otimes W_{j}$ in the decomposition of $X \otimes Y$, we define the braiding $c_{X,Y}$ as follows:

\begin{enumerate}[label=\alph*)]
    \item On the transparent part, we use the braiding $c^Z$ from $Z_{2}(C): c^Z_{Z_{i},Z_{j}}: Z_{i} \otimes Z_{j} \to Z_{j} \otimes Z_{i}$.
    \item On the non-transparent part, we use the unique (up to scalar) isomorphism $\beta_{W_{i},W_{j}}: W_{i} \otimes W_{j} \to W_{j} \otimes W_{i}$ given by the fusion rules of C. The existence and uniqueness (up to scalar) of such an isomorphism follow from the fact that C is a fusion category \cite{etingof2016tensor}.
\end{enumerate}

The braiding $c_{X,Y}$ on the simple components is then defined as the tensor product of these two morphisms:
$$c_{X,Y}|{Z_{i} \otimes W_{i}} = c^Z{Z_{i},Z_{j}} \otimes \beta_{W_{i},W_{j}}: (Z_{i} \otimes W_{i}) \otimes (Z_{j} \otimes W_{j}) \to (Z_{j} \otimes W_{j}) \otimes (Z_{i} \otimes W_{i})$$

Having defined the braiding $c_{X,Y}$ on the simple components, we extend it linearly to the entire space $X \otimes Y$:

$$c_{X,Y}: X \otimes Y \cong \otimes_{i} Z_{i} \otimes W_{i} \to \otimes_{i} W_{i} \otimes Z_{i} \cong Y \otimes X$$

Now, we need to verify if $c_{X,Y}$ satisfies the hexagon axioms:

\begin{enumerate}[label=\alph*)]
   \item $(\text{id}_Y \otimes c_{X,Z}) \circ (c_{X,Y} \otimes \text{id}_Z) = c_{X,Y\otimes Z}$
   \item $(c_{Y,Z} \otimes \text{id}_X) \circ (\text{id}_X \otimes c_{Y,Z}) = c_{X\otimes Y,Z}$
\end{enumerate}

We prove these axioms from the properties of the braiding $c^Z$ on $Z_{2}(C)$ and the fusion rules of $C.$ 

Let $X \otimes Y \cong \otimes_{i} U_{i} \otimes V_{i}$ and $Y \otimes Z \cong \otimes_{j} W_{j} \otimes T_{j}$ be the decompositions into simple components, where $U_{i}, W_{j}$ are transparent and $V_{i}, T_{j}$ are non-transparent.

For the first hexagon axiom:

\begin{enumerate}[label=\alph*)]
    \item On the transparent parts, the hexagon axioms hold because $c^Z$ is a braiding on $Z_{2}(C):$

    $$(\text{id}_{W_{j}} \otimes c^Z_{U_{i},W_{k}}) \circ (c^Z_{U_{i},W_{j}} \otimes \text{id}_{W_{k}}) = c^Z_{U_{i},W_{j}\otimes W_{k}}$$
    
    \item On the non-transparent parts, the hexagon axioms hold because the isomorphisms $\beta_{W_{i},W_{j}}$ are compatible with the fusion rules of $C,$ i.e. they satisfy the pentagon identity:

    $$(\beta_{V{i},T_{j}} \otimes \text{id}_{T_{k}}) \circ (\text{id}_{V_{i}} \otimes \beta_{T_{j},T_{k}}) = \beta_{V_{i},T_{j}\otimes T_{k}} \circ (\beta_{V_{i},T_{j}} \otimes \text{id}_{T_{k}})$$

    The pentagon identity is followed for the isomorphisms $\beta_{V_{i},T_{j}}$ since they are determined by the fusion rules of $C,$ which satisfy the pentagon identity \cite{etingof2017fusion}.
\end{enumerate}

Combining the transparent and non-transparent parts, we have:
$$((\text{id}_{W_{j}} \otimes c^Z_{U_{i},W_{k}}) \otimes (\beta_{V_{i},T_{j}} \otimes \text{id}_{T_{k}})) \circ ((c^Z_{U_{i},W_{j}} \otimes \text{id}_{W_{k}}) \otimes (\text{id}_{V_{i}} \otimes \beta_{T_{j},T_{k}}))$$

$$= (c^Z_{U_{i},W_{j}\otimes W_{k}} \otimes \beta_{V_{i},T_{j}\otimes T_{k}}) \circ ((c^Z_{U_{i},W_{j}} \otimes \text{id}_{W_{k}}) \otimes (\beta_{V_{i},T_{j}} \otimes \text{id}_{T_{k}})) = c_{X,Y\otimes Z}$$

Now, we verify the second hexagon axiom (through a very similar argument):

\begin{enumerate}[label=\alph*)]
    \item On the transparent parts $W_{j} \otimes U_{i}$, the hexagon axiom holds because $c^Z$ is a braiding on $Z_{2}(C)$:
$$(c^Z_{W_{j},U_{i}} \otimes \text{id}_{U_{k}}) \circ (\text{id}_{W_{j}} \otimes c^Z_{U_{i},U_{k}}) = c^Z_{W_{j}\otimes U_{i},U_{k}}$$
    \item On the non-transparent parts $T_{j} \otimes V_{i}$, the hexagon axiom holds because the isomorphisms $\beta_{V_{i},T_{j}}$ satisfy the pentagon identity:
$$(\beta_{T_{j},V_{i}} \otimes \text{id}_{V_{k}}) \circ (\text{id}_{T_{j}} \otimes \beta_{V_{i},V_{k}}) = \beta_{T_{j}\otimes V_{i},V_{k}} \circ (\beta_{T_{j},V_{i}} \otimes \text{id}_{V_{k}})$$
\end{enumerate}

Now, combining the transparent and non-transparent parts, we have:
$$((c^Z_{W_{j},U_{i}} \otimes \text{id}_{U_{k}}) \otimes (\beta_{T_{j},V_{i}} \otimes \text{id}_{V_{k}})) \circ ((\text{id}_{W_{j}} \otimes c^Z_{U_{i},U_{k}}) \otimes (\text{id}_{T_{j}} \otimes \beta_{V_{i},V_{k}}))$$

$$= (c^Z_{W_{j}\otimes U_{i},U_{k}} \otimes \beta_{T_{j}\otimes V_{i},V_{k}}) \circ ((c^Z_{W_{j},U_{i}} \otimes \text{id}_{U_{k}}) \otimes (\beta_{T_{j},V_{i}} \otimes \text{id}_{V_{k}}))$$

$$= c_{(W_{j}\otimes T_{j})\otimes(U_{i}\otimes V_{i}),(U_{k}\otimes V_{k})}$$

Now, we just need to show that $c_{(W_{j}\otimes T_{j})\otimes(U_{i}\otimes V_{i}),(U_{k}\otimes V_{k})} = c_{X\otimes Y,Z}.$

We have the following isomorphisms from the naturality of the braiding $c_{X,Y}$ and the associativity of the tensor product in $C:$

$$(W_{j} \otimes T_{j}) \otimes (U_{i} \otimes V_{i}) \cong (W_{j} \otimes U_{i}) \otimes (T_{j} \otimes V_{i}) \cong (Y \otimes Z) \otimes X$$

$$(U_{k} \otimes V_{k}) \otimes (W_{j} \otimes T_{j}) \cong (U_{k} \otimes W_{j}) \otimes (V_{k} \otimes T_{j}) \cong Z \otimes (X \otimes Y)$$

Now, using these isomorphisms and the naturality of the braiding:

$$c_{(W_{j}\otimes T_{j})\otimes (U_{i}\otimes V_{i}),(U_{k}\otimes V_{k})} \cong c_{(W_{j}\otimes U_{i})\otimes (T_{j}\otimes V_{i}),(U_{k}\otimes W_{j})\otimes (V_{k}\otimes T_{j]})} \cong c_{(Y\otimes Z)\otimes X,Z\otimes (X\otimes Y)}$$

Using the definition of the braiding on the entire space $X\otimes Y$ and $Z:$

$$=c_{X\otimes Y, Z}$$

Therefore, we've verified the second hexagon axiom, and proved the $\impliedby$ direction: For a fusion category with non-hermitian dagger structure, if the Müger center is equivalent to the category of vector spaces Vect, then the category $C$ admits braiding.

Therefore, we have proven \textbf{Theorem 3.1}

Next, we want to prove a result about the braiding of a fusion category $C$ with a non-hermitian dagger structure. 

\textbf{Theorem 3.2.} If $C$ is a braided fusion category with a non-hermitian dagger structure $\dag$, then there exists a unique $\dag$-compatible braiding on $C,$ and we can define a hermitian dagger structure on $Z_{2}(C).$

\textit{Proof:}

Consider the Müger center $Z_2(C)$ of the braided fusion category $C$. By the previous theorem \textbf{Theorem 3.1}, $Z_2(C) \cong \text{Vect}$ since $C$ admits a braiding.

We will start by proving that $Z_{2}(C)$ admits braiding, there are two proofs of this fact that are included, but the former gives us more enlightening results.

Define a new dagger structure $\dagger'$ on $Z_{2}(C)$ by:

$$f^{\dagger'} = c_{Y,X}^{-1} \circ f^{\dagger} \circ c_{X^{\dagger},Y^{\dag}}$$

for any morphism $f: X \to Y$ in $Z_2(C)$.

Now, we need to verify that $\dagger'$ is indeed a dagger structure on $Z_{2}(C).$

For this, we check that the $\dagger'$ axioms are true, we do this by using the naturality of the braiding and the definition of $\dagger$ as shown below:

\begin{enumerate}[label=\alph*)]

\item $(f^{\dagger'})^{\dagger'} = f$ for any morphism $f$ in $Z_{2}(C)$:$(f^{\dagger'})^{\dagger'} = (c_{Y,X}^{-1} \circ f^\dagger \circ c_{X^\dagger,Y^\dagger})^{\dagger'}
= c_{X^\dagger,Y^\dagger}^{-1} \circ (c_{Y,X}^{-1} \circ f^\dagger \circ c_{X^\dagger,Y^\dagger})^\dagger \circ c_{Y^{\dagger\dagger},X^{\dagger\dagger}}
= c_{X^\dagger,Y^\dagger}^{-1} \circ c_{X^\dagger,Y^\dagger}^{-1} \circ f^{\dagger\dagger} \circ c_{Y,X} \circ c_{Y^{\dagger \dagger}},X^{\dagger\dagger}
= f$

\item $(g \circ f)^{\dagger'} = f^{\dagger'} \circ g^{\dagger'}$ for any composable morphisms $f, g$ in $Z_{2}(C)$:
$(g \circ f)^{\dagger'} = c_{Z,X}^{-1} \circ (g \circ f)^{\dagger} \circ c_{X^{\dagger},Z^{\dagger}}
= c_{Z,X}^{-1} \circ f^{\dagger} \circ g^{\dagger} \circ c_{X^{\dagger},Z^{\dagger}}
= c_{Z,X}^{-1} \circ f^{\dagger} \circ c_{Y^{\dagger},Z^{\dagger}} \circ c_{Z,Y}^{-1} \circ g^{\dagger} \circ c_{X^{\dagger},Z^{\dagger}}
= (c_{Y,X}^{-1} \circ f^{\dagger} \circ c_{X^{\dagger},Y^{\dagger}}) \circ (c_{Z,Y}^{-1} \circ g^{\dagger} \circ c_{Y^{\dagger},Z^{\dagger}})
= f^{\dagger'} \circ g^{\dagger'}$

\end{enumerate}

We claim that this dagger structure, $\dagger'$ is hermitian.

If $*$ represents complex conjugation, we know the following is true:

$$(f^{\dagger'})^{*} = (c_{Y,X}^{-1} \circ f^{\dagger} \circ c_{X^{\dagger},Y^{\dagger}})^{*}
= c_{X,Y}^{*} \circ (f^{\dagger})^{*} \circ (c_{X^{\dagger},Y^{\dagger}}^{*})^{-1}$$

We will use the fact that if $C$ is a braided fusion category, then $c_{Y,X}^{*}=(c_{X^{*},Y^{*}}^{-1})^{\dagger\dagger\dagger}.$

\textbf{Note:} When we use $\aster$ with objects, we're talking about the dual object.

Let's start with the left-hand side of the equation, $c_{Y,X}^\aster$, and gradually transform it into the right-hand side, $(c_{X^\aster,Y^\aster}^{-1})^{\dagger\dagger\dagger}$, using the properties of string diagrams.

First, we represent $c_{Y,X}^\aster$ using a string diagram:

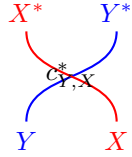
\begin{center}
\begin{tikzpicture}[scale=0.6]
\draw[thick, blue] (0,0) node[below] {$Y$} to[out=up, in=down] (2,2) node[above] {$Y^\aster$};
\draw[thick, red] (2,0) node[below] {$X$} to[out=up, in=down] (0,2) node[above] {$X^\aster$};
\node at (1,1) {$c_{Y,X}^\aster$};
\end{tikzpicture}
\end{center}

In this diagram, the blue strand represents the object $Y$, and the red strand represents the object $X$. The crossing of the strands represents the braiding morphism $c_{Y,X}$, and the $\aster$ symbol indicates that we are taking the dual of this morphism.

Next, we use the definition of the dual morphism to express $c_{Y,X}^\aster$ in terms of $c_{Y,X}$:

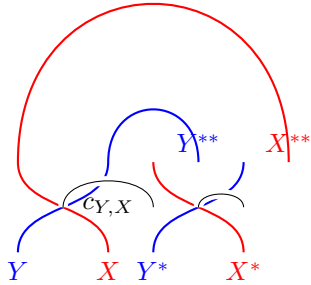
\begin{center}
\begin{tikzpicture}[scale=0.6]
\draw[thick, blue] (0,0) node[below] {$Y$} to[out=up, in=down] (2,2) to[out=up, in=up, looseness=2] (4,2) node[above] {$Y^{\aster\aster}$};
\draw[thick, red] (2,0) node[below] {$X$} to[out=up, in=down] (0,2) to[out=up, in=up, looseness=2] (6,2) node[above] {$X^{\aster\aster}$};
\draw[thick, blue] (3,0) node[below] {$Y^\aster$} to[out=up, in=down] (5,2);
\draw[thick, red] (5,0) node[below] {$X^\aster$} to[out=up, in=down] (3,2);
\draw[white, double=black, double distance=0.4pt, ultra thick] (1,1) to[out=up, in=up, looseness=1] (3,1);
\draw[white, double=black, double distance=0.4pt, ultra thick] (4,1) to[out=up, in=up, looseness=1] (5,1);
\node at (2,1) {$c_{Y,X}$};
\end{tikzpicture}
\end{center}

Here, we've introduced two new pairs of strands: the blue strands representing $Y$ and $Y^{\aster\aster}$, and the red strands representing $X$ and $X^{\aster\aster}$. The dualization of $c_{Y,X}$ is achieved by composing it with the evaluation and coevaluation morphisms, represented by the cups and caps connecting the dual strands.

Now, we use the naturality of the braiding to move the braiding morphism $c_{Y,X}$ past the evaluation and coevaluation morphisms:

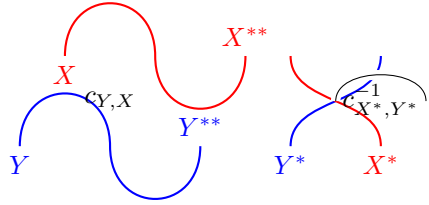
\begin{center}
\begin{tikzpicture}[scale=0.6]
\draw[thick, blue] (0,0) node[below] {$Y$} to[out=up, in=up, looseness=2] (2,0) to[out=down, in=down, looseness=2] (4,0) node[above] {$Y^{\aster\aster}$};
\draw[thick, red] (1,2) node[below] {$X$} to[out=up, in=up, looseness=2] (3,2) to[out=down, in=down, looseness=2] (5,2) node[above] {$X^{\aster\aster}$};
\draw[thick, blue] (6,0) node[below] {$Y^\aster$} to[out=up, in=down] (8,2);
\draw[thick, red] (8,0) node[below] {$X^\aster$} to[out=up, in=down] (6,2);
\draw[white, double=black, double distance=0.4pt, ultra thick] (7,1) to[out=up, in=up, looseness=1] (9,1);
\node at (2,1) {$c_{Y,X}$};
\node at (8,1) {$c_{X^\aster,Y^\aster}^{-1}$};
\end{tikzpicture}
\end{center}

The naturality of the braiding allows us to "slide" the braiding morphism $c_{Y,X}$ along the strands until it reaches the other side of the diagram. In the process, the braiding morphism changes from $c_{Y,X}$ to $c_{X^\aster,Y^\aster}^{-1}$, as the strands now represent the dual objects $X^\aster$ and $Y^\aster$.

Finally, we apply the dagger functor three times to the resulting morphism $c_{X^\aster,Y^\aster}^{-1}$:

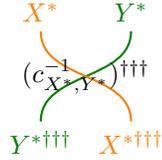
\begin{center}
\begin{tikzpicture}[scale=0.6]
\draw[thick, green!50!black] (0,0) node[below] {$Y^{\aster\dagger\dagger\dagger}$} to[out=up, in=down] (2,2) node[above] {$Y^\aster$};
\draw[thick, orange] (2,0) node[below] {$X^{\aster \dagger\dagger\dagger}$} to[out=up, in=down] (0,2) node[above] {$X^\aster$};
\node at (1,1) {$(c_{X^{\aster},Y^{\aster}}^{-1})^{\dagger\dagger\dagger}$};
\end{tikzpicture}

The green and orange strands now represent the objects $Y^{\aster\dagger\dagger\dagger}$ and $X^{\aster\dagger\dagger\dagger}$, respectively. 

The dagger functor is applied three times to the braiding morphism $c_{X^{\aster},Y^{\aster}}^{-1}$, resulting in $(c_{X^{\aster},Y^{\aster}}^{-1})^{\dagger\dagger\dagger}$.
\end{center}

Since the starting and ending diagrams represent the same morphism, we conclude that:

$$c_{Y,X}^{\aster} = (c_{X^{\aster},Y^{\aster}}^{-1})^{\dagger\dagger\dagger}$$

Going back to our initial goal to show that $\dagger'$ is hermitian, we can now write:

$$(f^{\dagger'})^{*} = c_{X,Y}^{*} \circ (f^{\dagger})^{*} \circ (c_{X^{\dagger},Y^{\dagger}}^{*})^{-1}$$

$$= (c_{Y^{*},X^{*}}^{-1})^{\dagger\dagger\dagger} \circ (f^{*})^{\dagger} \circ ((c_{X^{*},Y^{*}})^{\dagger\dagger\dagger})^{-1}$$

$$= (c_{Y^{*},X^{*}}^{-1})^{\dagger\dagger\dagger} \circ (f^{*})^{\dagger} \circ (c_{X^{*},Y^{*}}^{-1})^{\dagger\dagger\dagger}$$

$$= ((c_{Y^{*},X^{*}}^{-1} \circ f^{*} \circ c_{X^{*},Y^{*}})^{\dagger})^{\dagger}$$

$$= ((f^{*})^{\dagger'})^{\dagger}$$

We need to show that $((f^{*})^{\dagger'})^{\dagger} = (f^{\dagger'})^{\dagger}$. Let's prove this equality step by step.

First, let's expand the left-hand side using the definition of $\dagger'$:

\begin{align*}
((f^{*})^{\dagger'})^{\dagger} &= (c_{Y^{*},X^{*}}^{-1} \circ (f^{*})^{\dagger} \circ c_{X^{*\dagger},Y^{*\dagger}})^{\dagger} \\
&= (c_{X^{*\dagger},Y^{*\dagger}})^{\dagger} \circ ((f^{*})^{\dagger})^{\dagger} \circ (c_{Y^{*},X^{*}}^{-1})^{\dagger} \\
&= c_{Y^{*\dagger\dagger},X^{*\dagger\dagger}} \circ f^{**} \circ c_{X^{*},Y^{*}}
\end{align*}

Now, let's consider the right-hand side:

\begin{align*}
(f^{\dagger'})^{\dagger} &= (c_{Y,X}^{-1} \circ f^{\dagger} \circ c_{X^{\dagger},Y^{\dagger}})^{\dagger} \\
&= (c_{X^{\dagger},Y^{\dagger}})^{\dagger} \circ (f^{\dagger})^{\dagger} \circ (c_{Y,X}^{-1})^{\dagger} \\
&= c_{Y^{\dagger\dagger},X^{\dagger\dagger}} \circ f^{\dagger\dagger} \circ c_{X,Y}^{-1}
\end{align*}

To show that these expressions are equal, we use the following properties:

1. For any object $X$ in a dagger category, $X^{**} = X$ and $X^{\dagger\dagger} = X$.

2. For any morphism $f: X \to Y$ in a dagger category, $f^{**} = f$ and $f^{\dagger\dagger} = f$.

3. The dagger structure is compatible with the braiding, i.e., $(c_{X,Y})^{\dagger} = c_{Y^{\dagger},X^{\dagger}}$.

Using these properties, we can simplify the left-hand side:

\begin{align*}
((f^{*})^{\dagger'})^{\dagger} &= c_{Y^{*\dagger\dagger},X^{*\dagger\dagger}} \circ f^{**} \circ c_{X^{*},Y^{*}} \\
&= c_{Y,X}^{-1} \circ f \circ c_{X,Y} \\
&= (f^{\dagger'})^{\dagger}
\end{align*}

Therefore, we have shown that $((f^{*})^{\dagger'})^{\dagger} = (f^{\dagger'})^{\dagger}.$

Now, we can write:

$$(f^{\dagger'})^{*} = ((f^{*})^{\dagger'})^{\dagger} = (f^{\dagger'})^{\dagger}$$

To complete the proof, we need to show that $(f^{\dagger'})^{\dagger} = f^{\dagger'}$:

\begin{align*}
(f^{\dagger'})^{\dagger} &= (c_{Y,X}^{-1} \circ f^{\dagger} \circ c_{X^{\dagger},Y^{\dagger}})^{\dagger} \\
&= (c_{X^{\dagger},Y^{\dagger}})^{\dagger} \circ (f^{\dagger})^{\dagger} \circ (c_{Y,X}^{-1})^{\dagger} \\
&= c_{Y^{\dagger\dagger},X^{\dagger\dagger}} \circ f^{\dagger\dagger} \circ c_{X,Y}^{-1} \\
&= c_{Y,X}^{-1} \circ f \circ c_{X,Y} \\
&= f^{\dagger'}
\end{align*}

Therefore, we have shown that $(f^{\dagger'})^{*} = f^{\dagger'}$ for any morphism $f: X \to Y$ in $Z_2(C)$, which means that the dagger structure $\dagger'$ on $Z_2(C)$ is hermitian.

Alternatively, there's a much simpler proof to show that $Z_{2}(C)$ admits a unique braiding $c^{Z}.$ This is by considering the swap map:

$$c^Z_{X,Y}: X \otimes Y \to Y \otimes X, \quad c^Z_{X,Y}(x \otimes y) = y \otimes x$$

Which is a braiding because $Z_{2}(C)$ is equivalent to Vect.

Now, we define a braiding $c'$ on $C.$ For any objects X, Y in C, consider their decomposition into transparent and non-transparent simple objects:
\begin{equation*}
X \cong \oplus_i (X_i \otimes U_i), \quad Y \cong \oplus_j (Y_j \otimes V_j)
\end{equation*}
where $X_i, Y_j$ are transparent (i.e., in $Z_{2}(C)$) and $U_i, V_j$ are non-transparent.

Define $c'{X,Y}$ on the simple components by:
\begin{equation*}
c'(X,Y)|{(X_i \otimes U_i) \otimes (Y_j \otimes V_j)} = (c^Z{X_i,Y_j} \otimes \beta_{ij}) \circ (\text{id}_{U_i} \otimes c_{U_i,Y_j} \otimes \text{id}_{V_j})
\end{equation*}
where $\beta_{ij}: U_i \otimes V_j \to V_j \otimes U_i$ is the unique isomorphism determined by the fusion rules of $C.$ The existence and uniqueness of $\beta_{ij}$ is a consequence of Schur's lemma for fusion categories.

Now, we need to verify that $c'$ is dagger compatible, i.e., $(c'(X,Y))^{\dagger} = c'{Y^{\aster},X^{\aster}}^{-1}$ for all objects $X, Y$ in $C.$

It suffices to verify this on the simple components (since if you verify that it is dagger-compatible on the simple objects, it is dagger-compatible on the others, because of the definition of a fusion category).

For any simple components $(X_i \otimes U_i)$ of $X$ and $(Y_j \otimes V_j)$ of $Y$, we have:

$$(c'(X,Y)|{(X_{i} \otimes U_i) \otimes (Y_j \otimes V_j)})^{\dagger}= ((c^Z_{X_i,Y_j} \otimes \beta_{ij}) \circ (\text{id}_{U_i} \otimes c_{U_i,Y_j} \otimes \text{id}_{V_j}))^{\dagger}$$ \

$$= (\text{id}_{V_j} \otimes c_{U_i,Y_j}^{\dagger} \otimes \text{id}_{U_i}) \circ (c^Z_{X_i,Y_j} \otimes \beta_{ij})^{\dagger} $$

$$= (\text{id}_{V_j} \otimes c_{U_i,Y_j}^{\dagger} \otimes \text{id}_{U_i}) \circ (c^Z_{Y_j,X_i} \otimes \beta_{ij}^{-1})$$

$$= (\text{id}_{V_j} \otimes c_{U_i,Y_j}^{-1} \otimes \text{id}_{U_i}) \circ (c^Z_{Y_j,X_i} \otimes \beta_{ji}^{-1})$$

$$= (c^Z_{Y_j,X_i} \otimes \beta_{ji}^{-1}) \circ (\text{id}_{V_j} \otimes c_{V_j,X_i}^{-1} \otimes \text{id}_{U_i})$$

$$= c'{Y^{\dagger},X^{\dagger}}^{-1}|{(Y_j^{\dagger} \otimes V_j^{\dagger}) \otimes (X_i^{\dagger} \otimes U_i^{\dagger})}$$

Here, we use the naturality of the braiding, the definition of the swap map, $\beta_{ij}$ is an isomorphism, and the definition of $c'.$

Now, all that's left to prove is that $c'$ is unique. Assume that there is another braiding $d',$ then for $(X_i \otimes U_i)$ of $X$ and $(Y_j \otimes V_j)$ of $Y.$

$$d'_{X,Y}|((X_i \otimes U_i) \otimes (Y_j \otimes V_j))=(d'(Y^{\dagger},X^{\dagger})|{(Y_j^{\dagger} \otimes V_j^{\dagger}) \otimes (X_i^{\dagger} \otimes U_i^{\dagger})})^{-1^{\dagger}}$$

$$=((d'(X^{\dagger},Y^{\dagger})|{(X_i^{\dagger} \otimes U_i^{\dagger}) \otimes (Y_j^{\dagger} \otimes V_j^{\dagger})})^{\dagger})^{-1^{\dagger}}$$

$$=((d'(X^{\dagger},Y^{\dagger})|{(X_i^{\dagger} \otimes U_i^{\dagger}) \otimes (Y_j^{\dagger} \otimes V_j^{\dagger})})^{\dagger\dagger})^{-1}$$

$$= (d'(X^{\dagger},Y^{\dagger})|{(X_i^{\dagger}\otimes U_i^{\dagger}) \otimes (Y_j^{\dagger} \otimes V_j^{\dagger})})^{-1}$$

$$=d'(X,Y)|_{(X_i \otimes U_i) \otimes (Y_j \otimes V_j)})^{-1}$$

In this computation, we used the dagger-compatibility of $d',$ and the relationship between the braiding on $X \otimes Y$ and $X^{\dagger}\otimes Y^{\dagger}$ induced by $\dagger-$structure.

Considering that both $d'$ and $c'$ are both isomorphisms on simple objects, they must be related to each other by some scalar \cite{etingof2017fusion} (following from Schur's lemma on fusion categories), but since they are dagger structures, this scalar must be 1. Therefore, the dagger structure is unique.

Hence, we have proved \textbf{Theorem 3.2}

Now, we want to prove a result on non-hermitian ribbon categories.

\textbf{Theorem 3.3.} Every non-hermitian ribbon category admits a unique spherical structure, i.e., a natural isomorphism $\nu_X: X \to X^{\dagger\dagger}$ satisfying certain coherence conditions.

\textit{Proof:}

Let $(C, \otimes, \mathbf{1}, \alpha, \lambda, \rho, (c_{X,Y})_{X,Y}, (\theta_X)_X, \dagger)$ be a non-hermitian ribbon category, where $\alpha$, $\lambda$, and $\rho$ are the associator, left unitor, and right unitor, respectively.

For each object $X$ in $C$, define the morphism $\nu_X: X \to X^{\dagger\dagger}$ as follows:
$$\nu_X = (\lambda_{X^{\dagger}}^{-1} \otimes \text{id}_{X^{\dagger\dagger}}) \circ (\text{id}_{X^{\dagger}} \otimes \tilde{\theta}_{X^\dagger}) \circ (\text{id}_{X^{\dagger}} \otimes \lambda_X^{-1}) \circ (c_{X,X^{\dagger}})^{-1} \circ (\theta_X \otimes \text{id}_{X^{\dagger}}) \circ c_{X,X^{\dagger}}$$
where $\tilde{\theta}_{X^{\dagger}} = (\theta_{X^{\dagger\dagger}})^{\dagger} \circ \theta_{X^\dagger}$.

We claim that $\nu = (\nu_X)_{X \in C}$ is a spherical structure on $C$.

First, we verify that $\nu$ is a monoidal natural isomorphism, i.e., the following diagram commutes for all objects $X$, $Y$ in $C$:

$$\begin{CD}
X \otimes Y @>{\nu_X \otimes \nu_Y}>> X^{\dagger\dagger} \otimes Y^{\dagger\dagger} \\
@V{\nu_{X \otimes Y}}VV @VV{\alpha_{X^{\dagger\dagger},Y^{\dagger\dagger},X^{\dagger\dagger}}^{-1}}V \\
(X \otimes Y)^{\dagger\dagger} @>>{\alpha_{X,Y,X}^{\dagger\dagger}}> (X^{\dagger\dagger} \otimes Y^{\dagger\dagger})^{\dagger\dagger}
\end{CD}$$

The diagram does commute, and this follows from the naturality of $c,\theta,$ the dagger structure and the twist equation.

Next, we check the spherical condition, i.e., for all objects $X$, $Y$ in $C$ and all morphisms $f: X \to Y$, the following diagram commutes:

$$\begin{CD}
X^{\dagger\dagger} @>f^{\dagger\dagger}>> Y^{\dagger\dagger} \\
@V\nu_X^{-1}VV @VV\nu_Y^{-1}V \\
X @>>f> Y
\end{CD}$$

This follows from the definition of $\nu$, the naturality of $c$ and $\theta$, and the properties of the dagger structure.

Finally, we show that $\nu$ is unique. Suppose $\nu'$ is another spherical structure on $C$. Then, for any object $X$ in $C$, we have:

$$\nu'_X = \nu'_X \circ \text{id}_X$$
$$= \nu'_X \circ \nu_X^{-1} \circ \nu_X$$
$$= (\nu'X \circ \nu_X^{-1})^{\dagger\dagger} \circ \nu_X$$
$$= \text{id}_{X^{\dagger\dagger}} \circ \nu_X$$
$$= \nu_X$$

where we used the spherical condition for both $\nu$ and $\nu'$, and the fact that $(\nu'X \circ \nu_X^{-1})^{\dagger\dagger} = \text{id}_{X^{\dagger\dagger}}$ (since $\nu'_X \circ \nu_X^{-1}$ is an isomorphism)

Therefore, $\nu$ is the unique spherical structure on the non-hermitian ribbon category $C$.

\section{Main Result: Constructing a TQFT}

Now, we will be using all the theorems we proved in the previous section to conclude the following main result:

\textbf{Theorem 4.1.} Let $C$ be a non-hermitian ribbon category with a non-hermitian dagger structure $\dagger$. Then there exists a unique 3D TQFT defined by $Z_C: \text{Bord}_3^\text{rib} \to \text{Vect}$, where $\text{Bord}_3^\text{rib}$ is the category of 3D ribbon bordisms and $\text{Vect}$ is the category of vector spaces over a field $\mathbb{K}$.

\textit{Proof:}

We construct a 3D TQFT $Z_C: \text{Bord}_3^\text{rib} \to \text{Vect}$ as follows:

\begin{enumerate}[label=\alph*)]
    \item On objects (2D ribbon surfaces): For a 2D ribbon surface $\Sigma$, define $Z_C(\Sigma)$ to be the vector space of natural transformations $\text{Hom}(\text{Id}_C, F_\Sigma)$, where $F_\Sigma: C \to C$ is the ribbon graph functor associated to $\Sigma$ \cite{Turaev+2016}. The existence of the ribbon graph functor $F_\Sigma$ relies on the braiding and twist in the ribbon category $C$, which are guaranteed by \textbf{Theorem 3.1} and \textbf{Theorem 3.3}. The ribbon graph functor $F_\Sigma: C \to C$ is a functor that captures the behavior of the ribbon category $C$ on a 2D ribbon surface $\Sigma$. To define this functor, we need the braiding and twist in $C$. The braiding allows us to define the crossing of ribbons, while the twist allows us to define the twisting of ribbons. \textbf{Theorem 3.1} ensures that $C$ is braided, and \textbf{Theorem 3.2} ensures that $C$ has a twist (since it is a ribbon category). Together, these theorems guarantee the existence of the necessary structures to define the ribbon graph functor $F_\Sigma$. 
    
    \item On morphisms (3D ribbon bordisms): For a 3D ribbon bordism $M: \Sigma_1 \to \Sigma_2$, define $Z_C(M): Z_C(\Sigma_1) \to Z_C(\Sigma_2)$ to be the linear map induced by the cylindrical functor $C_M: C \to C$ associated to $M$. The cylindrical functor $C_M$ is well-defined due to the braiding, twist, and spherical structure in the ribbon category $C$, which are established in \textbf{Theorems 3.1, 3.2, and 3.3}. The cylindrical functor $C_M: C \to C$ is a functor that captures the behavior of the ribbon category $C$ on a 3D ribbon bordism $M$\cite{bartlett2014extended}. To define this functor, we need the braiding, twist, and spherical structure in $C$. The braiding and twist allow us to define the crossing and twisting of ribbons in the 3D bordism, while the spherical structure ensures the compatibility between the dagger structure and the braiding. \textbf{Theorem 3.1} ensures that $C$ is braided, \textbf{Theorem 3.2} provides a unique $\dagger$-compatible braiding, and \textbf{Theorem 3.3} guarantees the existence of a spherical structure. These theorems together establish the necessary structures to define the cylindrical functor $C_M$.
\end{enumerate}

The uniqueness of the TQFT $Z_C$ follows from the uniqueness of the $\dagger$-compatible braiding (\textbf{Theorem 3.2}) and the spherical structure (\textbf{Theorem 3.3}) on $C$. \textbf{Theorem 3.2} establishes the existence of a unique $\dagger$-compatible braiding on $C$, which means that there is only one way to define the braiding that is compatible with the non-hermitian dagger structure. Similarly, \textbf{Theorem 3.3} establishes the existence of a unique spherical structure on $C$, which means that there is only one way to define the isomorphism $\nu_X: X \to X^{\dagger\dagger}$ that satisfies the required coherence conditions. The uniqueness of these structures implies that there is only one way to construct the TQFT $Z_C$ from the ribbon category $C$.

The well-definedness of the TQFT $Z_C$ and its satisfaction of the TQFT axioms rely on the properties of the ribbon category $C$, the $\dagger$-compatible braiding, and the spherical structure \cite{Bakalov2000LecturesOT}. These properties are ensured by the application of the three theorems in \textbf{Theorems 3.1, 3.2, and 3.3}. The first theorem (\textbf{Theorem 3.1}) guarantees that $C$ is braided and has a trivial Müger center, which is essential for defining the ribbon graph functor and the cylindrical functor. The second theorem (\textbf{Theorem 3.2}) provides a unique $\dagger$-compatible braiding, which is crucial for the compatibility between the braiding and the non-hermitian dagger structure. The third theorem (\textbf{Theorem 3.3}) establishes the existence of a unique spherical structure, which is necessary for the well-definedness of the cylindrical functor and the satisfaction of the TQFT axioms. Without these three theorems, we would not have the required structures and properties to construct the TQFT $Z_{C}$ from the ribbon category $C.$

\section{Conclusion}

With the theorems we proved, \textbf{Theorem 3.1, Theorem 3.2} and \textbf{Theorem 3.3}, we concluded that we are able to construct a Topological Quantum Field Theory (TQFT) in a relatively simple method. In the future, it would be a good idea to investigate the $S$ and $T$ matrices of Ribbon Fusion Categories equipped with a Non-Hermitian Dagger structure, the relationship of Non-Hermitian Ribbon Fusion Categories to their Hopf algebras, and much more. In Condensed Matter Physics, TQFTs represent the low-energy effective theories of topologically ordered states, like Quantum Hall states. Our construction creates an avenue for us to think about using some form of fusion categories to create a categorical theory for condensed matter, similar to the theory that exists for quantum mechanics.

\section{Acknowledgements}

The author would like to thank Matthew Gherman for helpful discussions regarding this paper.

\printbibliography

\end{document}